\documentclass[a4paper]{article}
\usepackage{verbatim}
\usepackage{amsthm,amsmath,amssymb}
\usepackage{graphicx}
\usepackage{color}
\catcode`\@=11 \@addtoreset{equation}{section}

\catcode`\@=12

\newtheorem{Theorem}{Theorem}[section]
\newtheorem{Proposition}{Proposition}[section]
\newtheorem{Lemma}{Lemma}[section]
\newtheorem{Corollary}{Corollary}[section]
\newtheorem{Remark}{Remark}[section]
\newtheorem{Definition}{Definition}[section]

\newcommand{\bTheorem}[1]{
\begin{Theorem} \label{T#1} }
\newcommand{\eT}{\end{Theorem}}

\newcommand{\bProposition}[1]{
\begin{Proposition} \label{P#1}}
\newcommand{\eP}{\end{Proposition}}

\newcommand{\bLemma}[1]{
\begin{Lemma} \label{L#1} }
\newcommand{\eL}{\end{Lemma}}

\newcommand{\bCorollary}[1]{
\begin{Corollary} \label{C#1} }
\newcommand{\eC}{\end{Corollary}}

\newcommand{\bRemark}[1]{
\begin{Remark} \label{R#1} }
\newcommand{\eR}{\end{Remark}}

\newcommand{\bDefinition}[1]{
\begin{Definition} \label{D#1} }
\newcommand{\eD}{\end{Definition}}

\newcommand{\bFormula}[1]{
\begin{equation} \label{#1}}
\newcommand{\eF}{\end{equation}}

\newcommand{\Ov}[1]{\overline{#1}}

\newcommand{\DC}{C^\infty_c}

\newcommand{\vue}{\vu_\ep}

\newcommand{\RR}{\mathbb{R}}

\newcommand{\vu}{\vc{u}}
\newcommand{\vc}[1]{{\bf #1}}
\newcommand{\Div}{{\rm div}_x}
\newcommand{\Grad}{\nabla_x}

\newcommand{\dx}{{\rm d} {x}}
\newcommand{\dt}{{\rm d} t }

\newcommand{\intO}[1]{\int_{\Omega} #1 \ \dx}

\newcommand{\ep}{\varepsilon}

\definecolor{Cgrey}{rgb}{0.85,0.85,0.85}
\definecolor{Cblue}{rgb}{0.50,0.85,0.85}
\definecolor{Cred}{rgb}{1,0,0}
\definecolor{fancy}{rgb}{0.10,0.85,0.10}

\newcommand{\calM}{{\mathcal M}}
\newcommand{\calT}{{\mathcal T}}

\newcommand{\EEE}{\color{black}}

\newcommand{\GGG}{\color{black}}

\DeclareMathOperator{\deriv}{d}
\DeclareMathOperator{\loc}{loc}
\newcommand{\dit}{\deriv\!t}
\newcommand{\ddt}{\frac{\deriv\!{}}{\dit}}

\newcommand{\rhs}{right hand side}

\newcommand\Cbox[2]{%
    \newbox\contentbox%
    \newbox\bkgdbox%
    \setbox\contentbox\hbox to \hsize{%
        \vtop{
            \kern\columnsep
            \hbox to \hsize{%
                \kern\columnsep%
                \advance\hsize by -2\columnsep%
                \setlength{\textwidth}{\hsize}%
                \vbox{
                    \parskip=\baselineskip
                    \parindent=0bp
                    #2
                }%
                \kern\columnsep%
            }%
            \kern\columnsep%
        }%
    }%
    \setbox\bkgdbox\vbox{
        \color{#1}
        \hrule width  \wd\contentbox %
               height \ht\contentbox %
               depth  \dp\contentbox
        \color{black}
    }%
    \wd\bkgdbox=0bp%
    \vbox{\hbox to \hsize{\box\bkgdbox\box\contentbox}}%
    \vskip\baselineskip%
}

\definecolor{fuchsia}{rgb}{1.0, 0.0, 1.0}

\begin{document}


\title{Analysis of a diffuse interface model\\
of multispecies tumor growth}

\author{Mimi Dai\thanks{Department of Mathematics, University of Illinois at Chicago, 851 S.~Morgan Street,
Chicago, IL 60607-7045, USA. E-mail: \textit{mdai@uic.edu}.}
\and
Eduard Feireisl \thanks{Institute of Mathematics of the Academy of Sciences of the Czech Republic, 
\v Zitn\' a 25, CZ-115 67 Praha 1, Czech Republic. E-mail: \textit{feireisl@math.cas.cz}.}
\and
Elisabetta Rocca\thanks{Weierstrass Institute for Applied
Analysis and Stochastics, Mohrenstr.~39, D-10117 Berlin,
Germany. E-mail: \textit{rocca@wias-berlin.de} and Dipartimento di Matematica ``F. Enriques'',
Universit\`{a} degli Studi di Milano, Milano I-20133, Italy.
E-mail: \textit{elisabetta.rocca@unimi.it}.}
 \and
Giulio Schimperna\thanks{Dipartimento di Matematica ``F. Casorati'', Universit\`a degli Studi di Pavia,
via Ferrata 1, Pavia I-27100, Italy.
E-mail: \textit{giusch04@unipv.it}.}
\and Maria E. Schonbek\thanks{Department of Mathematics, University of California, 
Santa Cruz, CA 95064, USA. E-mail: \textit{schonbek@ucsc.edu}.}
}

\maketitle

\begin{abstract}
We consider  a diffuse interface model for tumor growth recently proposed in \cite{cwsl}. In this new approach
sharp interfaces are replaced by narrow transition layers arising due to adhesive forces among the cell species. Hence,
a continuum thermodynamically consistent model is introduced.  The resulting PDE system couples four different types of equations: 
a Cahn-Hilliard type equation for the tumor cells (which include proliferating and dead cells), a Darcy law for the tissue velocity
field, whose divergence may be different from 0 and depend on the other variables,
a transport equation for the proliferating (viable) tumor cells, and a  quasi-static reaction diffusion equation for the
nutrient concentration. We establish  existence of weak solutions for the PDE system coupled with suitable initial and boundary
conditions. In particular, the proliferation function at the boundary is supposed to be nonnegative on the set where the velocity $\vu$
satisfies $\vu\cdot\nu>0$, where $\nu$  is the outer  normal to the boundary of the domain.
We also study  a singular limit  as the diffuse interface coefficient tends to zero.

\end{abstract}

{\bf Key words:} tumor growth, diffuse interface model, Cahn-Hilliard equation, reaction-diffusion equation, Darcy law, existence of weak solutions, singular limits.
\vskip3mm
\noindent {\bf AMS (MOS) Subject Classification:} {35B25, 35D30, 35K35, 35K57, 35Q92, 74G25, 78A70, 92C17.}

\section{Introduction}
Mathematical modeling and analysis of tumor growth processes give 
important insights on cancer growth progression. The models are  expected to help 
to provide optimal treatment strategies. The behavior of tumors is a complex biological 
phenomenon, influenced by many factors, such as cell--cell and cell--matrix adhesion,
mechanical stress, cell motility and transport of oxygen, nutrients
and growth factors. In  recent  years, many mathematical
models of cancer have been proposed and various numerical simulations have 
been carried out (cf., e.g., the recent reviews \cite{byrne, cl, fasano, friedman}). 
A variety of models  are available to investigate different 
characteristics of cancer: single-phase continuum and
multiphase mixture models, and methods that combine both 
continuum and discrete components (cf., e.g., \cite[Chap.~7]{cl}).

We will address the problem of existence of weak solutions for 
a PDE system for a tumor growth model introduced in \cite{cwsl} 
(cf.\ also \cite{wlfc} and \cite{wlc}) and analyze a singular limit of
that model. The works listed above can be framed in the continuum tumor 
growth models category. This modeling approach has become central in the studies of tumor 
development in applied mathematics (cf.\ also \cite{am,oph}).
Actually, the translation of biological processes into models 
generally turns out to be simpler for discrete models than 
for continuum approaches. Nevertheless, discrete models can be
difficult to study analytically because the associated computational
cost rapidly increases with the number of cells modeled. This 
makes it difficult to simulate millimeter or greater sized tumors.
For this reason, in larger scale systems (millimeter to centimeter scale), continuum 
methods provide a good modeling alternative. Mixture models, 
on the other hand, provide the capability of simulating in 
detail the interactions among multiple cell species.

In the framework of continuum models, the diffuse interface method 
turns out to be particularly useful to describe multi-species tumor
growth processes. In this approach the sharp interfaces are replaced
by narrow transition layers arising due to the adhesion forces among
different cell-species. This choice is quite effective since it avoids
to introduce complicated boundary conditions across the tumor/host 
tissue and other species/species interfaces. This would have been 
the case  when considering sharp interface models. Moreover, the 
diffuse interface approach eliminates the need of tracking the 
position of the interfaces, which is one one of the main issues of such models.

The model derived in \cite{cwsl} consists of a Cahn-Hilliard system with transport and reaction terms which governs various types of cell concentrations. The reaction terms depend on the nutrient concentration (e.g., oxygen) which obeys to a quasi-static  advection-reaction-diffusion equation coupled to the Cahn-Hilliard equations. The cell velocities satisfy a generalized Darcy's  law where, besides the pressure gradient,  appears  also the so-called Korteweg force due to the cell concentration.

Numerical simulations of diffuse-interface models for tumor growth have 
been carried out in several papers (see, for instance, \cite[Chap.~8]{cl}
and references therein). However, a rigorous mathematical analysis of the 
resulting PDEs  is still in its beginning. To the best of our knowledge, 
the first related papers are concerned with a simplified model,
the so-called Cahn-Hilliard-Hele-Shaw system (see \cite{ltz}, cf.~also \cite{ww,wz}) 
in which the nutrient $n$, the source of tumor $S_T$ and the fraction 
$S_D$ of the dead cells are neglected. Moreover, very recent contributions
(see \cite{cgh, fgr, cgrs1, cgrs2}) are devoted to the analysis of a 
newly proposed  simpler model in \cite{hzo} (see also \cite{wzz}). In this 
model, velocities are set to zero and the state variables are reduced to the 
tumor cell fraction and the nutrient-rich extracellular water fraction.

In what follows we briefly introduce the model proposed in \cite{cwsl},  
where a complete description as well as numerical simulations are provided. 
Our multi-species tumor model includes the  mechanical interaction between
different species. The following notation will be used:
\begin{itemize}
\item $\phi_i, i=1,2, 3$:  {\em the volume fractions of  the cells: $\phi_1 =P$:  proliferating  cell fraction; $\phi_2= 
\phi_D$: dead cell fraction; $\phi_3 = \phi_H$: host cell fraction;} 
\item   $\Pi$: {\em  the cell-to-cell pressure};
\item   $\vu$:=$\vu_i, i=1,2,3$: {\em the tissue velocity field.  We assume  that the cells are tightly packed and they march together};
\item  $n$:  {\em the nutrient concentration};
\item $\Phi = \phi_D +P $: {\em the  volume fraction of the tumor cells 
which is split into the sum of the dead tumor cells 
and of the proliferating cells};
\item ${\bf J}_i$: {\em the fluxes that account for mechanical interactions among the species};
\item  $S_i, i=1,2,3$: {\em  account for inter-component mass exchange
as well as gains due to proliferation of cells and loss due to cell death.}
\end{itemize}
The variables above  are naturally constrained by the relation $\phi_H+\Phi=1$.

\noindent The volume fractions obey the mass conservation (advection-reaction-diffusion) equations:
\bFormula{i1}
\partial_t\phi_i+\Div(\vu\phi_i)=-\Div{\bf J}_i+\Phi S_i.
\eF
We have assumed that the densities of the components are matched.
Notice  that unlike in  \cite{cwsl}, for simplicity, the variable $\phi_W$ 
standing for the volume fraction of water has been omitted.
The total energy adhesion,  supposed  independent of $\phi_H$, has the form 
\[
E=\intO{\left(\mathcal{F}(\Phi)+\frac{1}{2}|\Grad\Phi|^2\right)},
\]
where $\mathcal{F}$ is a logarithmic type mixing potential (cf.~(\ref{ee8}) 
in Subsection~\ref{ass}).
Then, we define the fluxes ${\bf J}_\Phi$  and ${\bf J}_H$ as follows: 
\begin{align}\nonumber
&
{\bf J}_\Phi={\bf J}_1+{\bf J}_2:=-\Grad\left(\frac{\delta E}{\delta \Phi}\right)=-\Grad\left(\mathcal{F}'(\Phi)-\Delta\Phi\right):=-\Grad\mu,\\
\nonumber
&{\bf J}_H={\bf J}_3:=-\Grad\left(\frac{\delta E}{\delta \phi_H}\right)=\Grad\left(\frac{\delta E}{\delta \Phi}\right),
\end{align}
where we have used in the last equality the fact that $\phi_H=1-\Phi$ and where
$\mu$ is the chemical potential of the system. For the source of mass 
in the host tissue  we have the following relations:
\begin{itemize}
\item $S_T=S_D+S_P:=S_2+S_1$,
\item $\Phi S_H:=\Phi S_3=\phi_H S_T=(1-\Phi)S_T$.
\end{itemize}
Assuming the mobility of the system to be constant, then the
tumor volume fraction $\Phi$ and the host tissue volume fraction~$\phi_H$
obey the following mass conservation equations (cf.~(\ref{i1})):
\begin{align}\label{massfi}
&\partial_t \Phi + \Div(\vu  \Phi) =-\Div{\bf J}_\Phi+\Phi(S_2+S_1),\\
\label{massfih}
&\partial_t \phi_H + \Div(\vu  \phi_H) =-\Div{\bf J}_H +\Phi S_3.
\end{align}
Using now the fact that $S_T=S_1+S_2$ and recalling that $\phi_H+\Phi=1$, 
we can forget of the equation for $\phi_H$ and we recover the equation 
for $\Phi$ in the form 
\begin{equation}\label{e:Fi}
\partial_t \Phi + \Div(\vu  \Phi) -\Div (\Grad \mu ) =  \Phi S_T , \ \mu =  \mathcal{F}'(\Phi)- \Delta \Phi.
\end{equation}
As  in \cite{wlfc}, we suppose the net source of tumor cells $S_T$ 
to be given by
\[
S_T=S_T(n, P, \Phi)=\lambda_MnP-\lambda_L(\Phi-P),
\]
where $\lambda_M\geq 0$ is the mitotic rate 
and $\lambda_L\geq 0$ is the lysing rate of dead cells. The volume 
fraction of dead tumor cells $\phi_D$ would satisfy an equation
similar to~\eqref{e:Fi}, namely
\[
\partial_t \phi_D + \Div(\vu  \phi_D) - \Div (\Grad \mu ) =  \Phi S_D,
\]
where the source of dead cells is taken as
\[
S_D=S_D (n,P,\Phi) = \left( \lambda_A + \lambda_N H (n_N - n) \right) P - \lambda_L (\Phi - P ).
\]
However, we prefer to couple the equation for $\Phi$ with the one for 
$P=\Phi-\phi_D$ which then reads
\[
\partial_t P + \Div (\vu P) = \Phi (S_T - S_D).
\]
Here $\lambda_AP$ describes the death of cells due to apoptosis 
(cf.~\cite[p.~730]{cwsl}) with rate $\lambda_A\geq 0$ and 
the term $\lambda_NH(n_N-n)P$ models the death of cells due 
to necrosis with rate $\lambda_N\geq 0$. In \cite{wlfc}
$H$ was originally taken as the Heaviside function. Here, for mathematical
reasons, we smooth it out by taking it as a regular and nonnegative
function of~$n$.
The term $n_N$ represents the necrotic limit, at which the tumor
tissue dies due to lack of nutrients.

The tumor velocity field $\vu$ (given by the mass-averaged velocity 
of all the components) is assumed to fulfill Darcy's law:
\[
\vu = - \Grad \Pi + \mu \Grad \Phi,
\]
where,  for simplicity, the motility has been taken  constant and equal to 1.
Summing up equations (\ref{i1}), we end up with the following constraint for the velocity field:
\[
\Div\vu=S_T.
\]
Since the time scale for nutrient diffusion is much faster than the rate of 
cell proliferation, the nutrient is assumed to evolve quasi-statically:
\[
- \Delta n + \nu_UnP = T_c (n, \Phi),
\]
where the nutrient capillarity term $T_c$ is
\[
T_c (n, \Phi) = \left[ \nu_1 (1 - Q(\Phi)) + \nu_2 Q(\Phi) \right] (n_c - n),
\]
$\nu_U$ represents the nutrient uptake rate by the viable tumor cells,
$\nu_1,$ $\nu_2$ denote the nutrient transfer rates for preexisting vascularization in the tumor
and host domains, and $n_c$ is the nutrient level of capillaries. 
The function $Q(\Phi)$ is assumed to be regular and to satisfy
$\nu_1 (1 - Q(\Phi)) + \nu_2 Q(\Phi) \ge 0$
(cf.~(\ref{ee13-}) below).
\begin{Remark} 
We chose the boundary conditions  proposed in \cite{cwsl} for $\Phi$, $\mu$, $\Pi$ and $n$.
On the other hand, under the homogeneous Neumann boundary conditions suggested in \cite{cwsl}
for $P$, we could not show that the system is well-posed. 
For this reason, we chose the boundary conditions (\ref{e7}), which 
are natural in connection with the transport equation (\ref{e4}) for $P$.
In particular, the proliferation function at the boundary has to 
be nonnegative on the set where the velocity $\vu$ satisfies $\vu\cdot\nu>0$, 
with $\nu$ denoting the outer normal unit vector to the boundary of our domain~$\Omega$.
By maximum principle, this implies in particular that $P\geq 0$ in $\Omega$,
which is an information we need for proving well-posedness of the system.
\end{Remark}
\noindent%
In summary, let $\Omega\subset\RR^3$ be a bounded domain and $T>0$ 
the final time of the process. For simplicity, choose $\lambda_M=\nu_U=1$, 
$\lambda_A=\lambda_1$, $\lambda_N=\lambda_2$, $\lambda_L=\lambda_3$. Then, 
in $\Omega\times (0,T)$, we have the following  system of equations:
\bFormula{e1}
\partial_t \Phi + \Div(\vu  \Phi) - \Div (\Grad \mu ) =  \Phi S_T , \ \mu = - \Delta \Phi + \mathcal{F}'(\Phi),
\eF
\bFormula{e2}
\vu = - \Grad \Pi + \mu \Grad \Phi,
\eF
\bFormula{e3}
\Div \vu = S_T,
\eF
\bFormula{e4}
\partial_t P + \Div (\vu P) = \Phi (S_T - S_D),
\eF
\bFormula{e5}
- \Delta n + nP = T_c (n, \Phi),
\eF
where
\bFormula{def:ST}
S_T (n,P,\Phi) = nP - \lambda_3 (\Phi - P),
\eF
\[
S_D (n,P,\Phi) = \left( \lambda_1 + \lambda_2 H (n_N - n) \right) P - \lambda_3 (\Phi - P ),
\]
\bFormula{def:Tc}
T_c (n, \Phi) = \left[ \nu_1 (1 - Q(\Phi)) + \nu_2 Q(\Phi) \right] (n_c - n).
\eF
The functions $Q, H$ and the constants $\lambda_i, \nu_i$ 
will be described in Section~\ref{ass}. 
System (\ref{e1}--\ref{e5}) will be coupled  with the 
following boundary conditions on $\partial\Omega\times (0,T)$:
\bFormula{e6}
\mu = \Pi = 0 , \ n = 1,
\eF
\bFormula{e7}
\Grad \Phi \cdot \nu = 0, \ \GGG P \vu \cdot \nu \geq 0, \EEE
\eF
and with the initial conditions
\bFormula{e8}
\Phi(0)=\Phi_0, \quad P(0)=P_0\ \mbox{in}\ \Omega.
\eF
Note that, as $P \geq 0$, the second condition in (\ref{e7}) 
should be interpreted as $P = 0$ whenever $\vu \cdot \nu < 0$, 
meaning on the part of the inflow part of the boundary.
\GGG Moreover, in the weak formulation, that condition 
will be incorporated into equation \eqref{e4} turning it
into a variational inequality (cf.~\eqref{wf3} below). \EEE

The different nature of the four equations as well as their nonlinear coupling
(especially due to the Korteweg term in the pressure equation) make the analysis 
of the problem particularly challenging. Moreover, we may notice that the singular
limit studied in the last Section~\ref{S} as the interface energy coefficient
is let tend to zero can be obtained only under more restrictive assumptions 
on the potential ${\mathcal{F}}$, which is required to be strictly convex,
\GGG and under different boundary conditions for $\vu$ (namely, 
we assume no-flux, rather than Dirichlet, conditions for $\Pi$). \EEE
We refer the reader to Remark~\ref{RR1} below for further comments
and for the discussion of related open problems.

\medskip

\noindent
{\bf Plan of the paper.} The main results and assumptions are stated in Section~\ref{M}.
The subsequent Sections~\ref{B} and \ref{SS} are the core of the paper where we provide
the a priori bounds for our solutions and we show the weak sequential stability properties.
In Section~\ref{A}, we construct an approximation scheme compatible with the apriori
estimates and prove its well-posedness. In the last section, we analyze the singular
limit problem mentioned above.


\section{Assumptions and main results}
\label{M}


\subsection{Singular potential and initial data}
\label{ass}

We suppose that the potential $\mathcal{F}$ supports the natural bounds 
\bFormula{ee7}
0 \leq \Phi(t,x) \leq 1.
\eF
To this end, we take $\mathcal{F} = \mathcal{C} + \mathcal{B}$, where $\mathcal{B} \in C^2(\RR)$ and
\bFormula{ee7b}
\mathcal{C} : \RR \mapsto [0, \infty] \ \mbox{convex, lower-semi continuous}, \ \mathcal{C}(\Phi) = \infty
\ \mbox{for}\ \Phi < 0 \ \mbox{or}\ \Phi > 1. 
\eF
\GGG Moreover, we ask that
\bFormula{ee7c}
 \mathcal{C} \in C^1(0,1), \ 
  \lim_{\Phi\to 0^+} \mathcal{C'}(\Phi) = \lim_{\Phi\to 1^-} \mathcal{C'}(\Phi) = \infty.
\eF  \EEE
A typical example of such $\mathcal{C}$ is the {\em logarithmic potential}
\bFormula{ee8}
\mathcal{C}(\Phi) = \left\{ \begin{array}{l} \Phi \log (\Phi) + (1 - \Phi) \log( 1 - \Phi) \ \mbox{for}\ \Phi \in [0,1],\\ \\
\infty \ \mbox{otherwise.} \end{array} \right.
\eF
\GGG \bRemark{Rsing}
 Condition~\eqref{ee7c} has mainly a technical character and is assumed just for the purpose
 of constructing a not too complicated approximation scheme (cf.~also Remark~\ref{Rrem:sep}). At the price of some additional 
 technical work it could be avoided. One may, for instance, consider the case where
 $\mathcal{C}(\Phi) = I_{[0,1]}(\Phi)$ (the {\em indicator function}\/ of $[0,1]$), which
 does not satisfy~\eqref{ee7c}.
\eR \EEE
\noindent%
Regarding the functions $Q$ and $H$ and the constants $\lambda_i$, $\nu_i$ appearing in the definitions of $S_T$ and $S_D$, we 
assume \GGG $Q, H \in C^1(\RR)$ together with \EEE
\bFormula{aHl}
\lambda_i \geq 0  \ \mbox{for}\ i=1,2,3, \ \ H \geq 0.
\eF
\bFormula{aQ}
\left[ \nu_1 (1 - Q(\Phi)) + \nu_2 Q(\Phi) \right] \geq 0, \ \ 0 < n_c < 1.
\eF
Finally, we suppose $\Omega$ be a bounded domain with smooth boundary in $\RR^3$ and impose the following conditions on the initial data: 
\bFormula{iphi}
\Phi_0\in H^1(\Omega), \quad 0 \leq \Phi_0 \leq 1, \quad \mathcal{C}(\Phi_0)\in L^1(\Omega),
\eF
\bFormula{iP}
P_0\in L^2(\Omega), \quad 0 \leq P_0 \leq 1\, \quad\hbox{a.e. in }\Omega.
\eF


\subsection{Main result}

Before stating the main result, let us introduce a suitable weak formulation of the problem. We say that 
$(\Phi, \vu, P, n)$ is a weak solution to problem (\ref{e1}--\ref{e8}) in $(0,T) \times \Omega$ if 
\begin{itemize}
\item[\bf (i)] these functions 
belong to the regularity class:
\bFormula{rid1}
\Phi \in C^0([0,T];H^1(\Omega))\cap L^2(0,T; W^{2,6}(\Omega)),
\eF
\bFormula{rid1b}
 \mathcal{C}(\Phi) \in L^\infty(0,T; L^1(\Omega)),
  \ \text{ hence, in particular,}\ 0 \leq \Phi \leq 1 \ \mbox{a.a. in} \ (0,T) \times \Omega;
\eF
\bFormula{rid2}
 \vu\in L^2((0,T)\times\Omega; \RR^3), \ {\rm div} \ \vu\in L^\infty((0,T)\times\Omega);
\eF
\bFormula{rid2bis}
\Pi\in L^2(0,T;W^{1,2}_0(\Omega)), \quad \mu\in \GGG L^2 \EEE(0,T; W^{1,2}_0(\Omega));
\eF
\bFormula{rid3}
P\in L^\infty((0,T)\times\Omega), \, 0 \leq P \leq 1 \  \mbox{ a.a. in} \ (0,T) \times \Omega;
\eF
\bFormula{rid4}
 n\in L^2(0,T;W^{2,2}(\Omega)), \ 0 \leq n \leq 1 \ \mbox{ a.a. in} \ (0,T) \times \Omega;
\eF
\item[\bf (ii)] the following integral identities hold:
\bFormula{wf1}
 \int_0^T \intO{ \left[ \Phi \partial_t \varphi + \Phi \vu \cdot \Grad \varphi + \mu \Delta \varphi + \Phi S_T \varphi \right] } \ \dt
  = - \intO{ \Phi_0 \varphi(0,\cdot) } 
\eF
for any $\varphi \in \DC([0,T) \times \Omega)$, where
\bFormula{wf2}
\mu = - \Delta \Phi + \mathcal{F}'(\Phi),\ \vu = - \Grad \Pi + \mu \Grad \Phi,
\eF
\bFormula{wf2bis} 
\Div \vu = S_T\ \mbox{a.a. in} \ (0,T) \times \Omega; \quad \Grad\Phi\cdot\nu|_{\partial\Omega}=0;
\eF  
\bFormula{wf3}
\int_0^T \intO{ \left[ P \partial_t \varphi + P \vu \cdot \Grad \varphi + \Phi (S_T - S_D ) \varphi \right] } \ \dt \geq - \intO{ P_0 \varphi(0, \cdot) }
\eF
for any $\varphi \in \DC([0,T) \times \Ov{\Omega})$, $\varphi|_{\partial \Omega} \geq 0$;
\bFormula{wf4}
-\Delta n + nP = T_c(n, \Phi) \ \mbox{a.a. in}\ (0,T) \times \Omega;\ n |_{\partial \Omega} = 1.
\eF
\end{itemize}
Now, we are able to state the main result of the present paper:
\bTheorem{theoremmain}
Let $T > 0$ be given. Under the assumptions stated in Subsection~\ref{ass}, the variational 
formulation (\ref{wf1}--\ref{wf4}) of the initial-boundary value problem 
(\ref{e1}--\ref{e8}) admits at least one solution in the regularity class 
(\ref{rid1}--\ref{rid4}). 
\eT
\bRemark{rembound}
 It is worth observing once more that the second boundary condition~\eqref{e7} is
 now incorporated into the variational inequality~\eqref{wf3}.
\eR


\section{A priori bounds}
\label{B}

In this section we establish several formal a priori estimates for our solution. The procedure
turns out to be rigorous when (smoother) solutions of the
approximated problem (\ref{e1r}--\ref{e5r}) are considered. 
In particular, this happens for the regularized solution constructed in 
Section~\ref{A} below. In this section we refer to system (\ref{e1}--\ref{e8}) and not to the 
weak formulation (\ref{wf1}--\ref{wf4}) because actually the a-priori estimates should be performed on 
the regularized problem (\ref{e1r}--\ref{e5r}) whose solutions are more regular than the ones obtained at the limit. 

We start with noticing that, as a direct consequence of our 
choice of the potential $\mathcal{F}$, the phase field function $\Phi$ 
satisfies (\ref{ee7}).


\subsection{Lower bound for $P$}

The density function $P$ satisfies the transport equation (\ref{e4}), which can be equivalently
rewritten in the form
\bFormula{ee9}
\partial_t P + \vu \cdot \Grad P = - P S_T + \Phi (S_T - S_D)
\eF
\[
= P \left[ - S_T + \Phi \left( n - \left( \lambda_1 + \lambda_2 H (n_N - n) \right) \right)
\right].
\]
Thus, provided 
\[ 
P(0, \cdot) = P_0 \geq 0,\; \mbox{and} \;P (t,x) \geq 0 \ \mbox{for}\ x \in \partial \Omega, \ \vu \cdot \nu \leq 0,
\]
we can deduce by maximum principle arguments that
\[
P \geq 0.
\]


\subsection{Positivity and upper bound for $n$}

In order to obtain positivity of $n$ we need
\[
- nP + T_c (n, \varphi) = -nP + \left[ \nu_1 (1 - Q(\Phi)) + \nu_2 Q(\Phi) \right] (n_c - n)
\]
to be positive (non-negative) whenever $n < 0$; actually, this follows from the hypothesis
(cf.~\eqref{aQ} in Subsection~\ref{ass})
\bFormula{ee13-}
\left[ \nu_1 (1 - Q(\Phi)) + \nu_2 Q(\Phi) \right] \geq 0, \ 0 < n_c < 1.
\eF
This assumption also implies that $n \leq 1$, so we may conclude that
\bFormula{ee13--}
 0 \leq n(t,x) \leq 1.
\eF


\subsection{Upper bound for $P$}

Since $0\leq \Phi\leq 1$ and $0\leq n\leq 1$, by the assumptions provided in~Subsection~\ref{ass}
we have 
\[
- \Phi \left( \lambda_1 + \lambda_2 H (n_N - n) \right) \leq 0.
\]
Hence evaluating the expression on the right-hand side of (\ref{ee9}) for $P = 1$ yields 
\[
 P \left[ - S_T + \Phi \left( n - \left( \lambda_1 + \lambda_2 H (n_N - n) \right) \right)\right] \leq
  \lambda_3 (\Phi - 1) + n (\Phi - 1).
\]
Consequently, provided  
\[ 
0 \leq P(0, \cdot) = P_0 \leq 1,\,\mbox{and }0 \leq P (t,x) \leq 1 \ \mbox{for}\ x \in \partial \Omega, \ \vu \cdot \nu \leq 0,
\]
it follows that
\bFormula{ee13}
0 \leq P(t,x) \leq 1.
\eF


\subsection{Estimates for the Cahn-Hilliard equation}
\label{ECH}

The standard estimates are obtained via multiplication of (\ref{e1}) by $\mu$:
\bFormula{ee14}
\frac{{\rm d}}{{\rm d}t} \intO{ \left[ \frac{1}{2} |\Grad \Phi |^2 + \mathcal{F}(\Phi) \right] } +
\intO{ |\Grad \mu|^2 } = - \intO{ \vu \cdot \Grad \Phi \mu },
\eF
where, by virtue of (\ref{e2}),
\[
- \intO{ \vu \cdot \Grad \Phi \mu } = - \intO{ |\vu|^2 } + \intO{ \Pi \Div \vu }
= - \intO{ |\vu|^2 } + \intO{ \Pi S_T }.
\]
Consequently, (\ref{ee14}) reads
\bFormula{ee15}
\frac{{\rm d}}{{\rm d}t} \intO{ \left[ \frac{1}{2} |\Grad \Phi |^2 + \mathcal{F}(\Phi) \right] } +
\intO{ \left[ |\Grad \mu|^2 + |\vu|^2 \right] } =  \intO{ \Pi S_T },
\eF
where
\[
  \left| \intO{ \Pi S_T } \right| 
    \leq \| S_T \|_{L^\infty(\Omega)} \| \Pi \|_{L^1(\Omega)}.
\]
Seeing that $\Pi$ solves the Dirichlet problem
\[
- \Delta \Pi = S_T - \Div (\mu \Grad \Phi), \ \Pi|_{\partial \Omega} = 0,
\]
we deduce that
\[
\| \Pi (t, \cdot) \|_{H^1(\Omega)} \leq \| S_T (t, \cdot) \|_{L^2(\Omega)} + \| \mu \Grad \Phi \|_{L^2(\Omega; \RR^3)},
\]
where, by means of Gagliardo-Nirenberg interpolation inequality,
\[
\| \mu \Grad \Phi \|_{L^2(\Omega; \RR^3)} \leq \| \mu (t, \cdot) \|_{L^4 (\Omega)} \| \Grad \Phi \|_{L^4(\Omega;\RR^3)}
\]
\[
\leq c \| \mu (t, \cdot) \|_{L^4 (\Omega)} \| \Phi (t, \cdot) \|_{L^\infty(\Omega)}^{1/2} \| \Delta \Phi (t, \cdot) \|_{L^2(\Omega)}^{1/2}
\]
\[
\leq c \| \mu (t, \cdot) \|_{L^4 (\Omega)} \| \Phi (t, \cdot) \|_{L^\infty(\Omega)}^{1/2}
  \left( \| \mu \|_{L^2(\Omega)}^{1/2}+\|\nabla\Phi\|_{L^2(\Omega)}^{1/2}\right),
\]
where the last inequality has been obtained testing the second~\eqref{e1} by~$\Phi$
and using the properties of~$\mathcal{F}$ (in particular, the monotonicity of $\mathcal{C}'$).

Thus, going back to (\ref{ee15}) and applying a standard version of Gr\"onwall's lemma, we deduce the bounds
\bFormula{ee16}
\sup_{t \in (0,T)} \| \Phi \|_{H^1(\Omega)} \leq c,
\eF
\bFormula{ee17}
\int_0^T \left[ \| \Grad \mu \|^2_{L^2(\Omega;\RR^3)} + |\vu|^2 \right] \ \dt \leq c.
\eF


\subsubsection{More estimates on $\Phi$}

Knowing that
\bFormula{ee18}
- \Delta \Phi + \mathcal{C}'(\Phi) = g = \mu - \mathcal B'(\Phi)  \in L^2(0,T; H^1(\Omega)),
\eF
we may multiply this relation by $- \Delta \Phi$ and use once more the 
monotonicity of $\mathcal C'$ to deduce
\[
\int_0^T \|\Phi \|^2_{W^{2,2}(\Omega)} \ \dt \leq c.
\]
Next, take an increasing function $h$ and multiply (\ref{ee18}) by $h (\mathcal{C}'(\Phi)) $ to obtain
\bFormula{ee181}
   \intO{ \left[ h'(\mathcal{C}'(\Phi)) \mathcal{C}''(\Phi) |\Grad \Phi |^2 
    +  h (\mathcal{C}'(\Phi)) \mathcal{C}'(\Phi) \right] } = \intO{ g h (\mathcal{C}'(\Phi)) }.
\eF
Choosing $h(\cdot) = (\cdot)^5$ and using that
$g \in L^2(0,T; L^6(\Omega))$, we then easily deduce
\[
\mathcal{C}'(\Phi) \ \mbox{is bounded in}\ L^2(0,T; L^6(\Omega)),
\]
whence, comparing terms in \eqref{ee18}, we also infer
\bFormula{ee19}
 \int_0^T \| \Phi \|^2_{W^{2,6}(\Omega)} \ \dt \leq c.
\eF


\subsubsection{Estimates on $\vu$}

Note that we already know
\[
\Div \vu = S_T \ \mbox{bounded in}\ L^\infty((0,T) \times \Omega)
\]
and
\[
\vu \ \mbox{bounded in}\ L^2((0,T) \times \Omega; \RR^3).
\]
Next, we compute
\[
 \GGG {\bf curl}_x \vu = \Grad \mu \wedge \Grad \Phi \EEE
  \in L^2(0,T; L^1(\Omega)) \cap L^1(0,T; L^2(\Omega)).
\]
\GGG Hence, we may take a test function $\varphi\in C^\infty(\RR^3)$ 
with support contained in $\Omega$ and apply~\cite[p.~51]{gr} to the 
function $\varphi\vu$. In view of the fact that
$\Div(\varphi\vu)$ and ${\bf curl}(\varphi\vu)$ are bounded in $L^1(0,T;L^2(\RR^3))$, 
we then obtain that $\varphi \vu$ is bounded in $L^1(0,T;H^1(\RR^3))$.
Consequently, $\vu$ satisfies
\bFormula{ee20}
  \int_0^T \| \vu \|_{H^1_{\loc}(\Omega; \RR^3)} \ \dt
\eF
\EEE


\section{Weak sequential stability}

\label{SS}

Suppose that
\[
\{ \Phi_\delta, \vu_\delta, P_\delta, n_\delta \}_{\delta > 0}
\]
is a family of solutions complying with the {\em a priori}\/ bounds obtained in the last section.  
Our goal is to show the precompactness of this family of solutions, that is to prove that
\[
\left\{ \begin{array}{c}
\Phi_\delta \to \Phi \ \mbox{weakly-(*) in} \ L^\infty((0,T) \times \Omega),\\
\ \vu_\delta \to \vu \ \mbox{weakly in}\ L^2((0,T) \times \Omega; \RR^3),\\ \ P_\delta \to P \ \mbox{weakly-(*) in}\ L^\infty((0,T) \times \Omega) ,\\ \ n_\delta \to n
\ \mbox{weakly-(*) in}\ L^\infty((0,T) \times \Omega),
\end{array}
\right\}
\]
where the limits solve the same system of equations.


\subsection{Compactness of the time derivatives}

It follows from (\ref{e1}) and the a-priori estimates we have on $\Phi$ that
\[
\partial_t \Phi_\delta \to \partial_t \Phi \ \mbox{weakly in}\ L^2(0,T; W^{-1,2}(\Omega)),
\]
whence, in accordance with (\ref{ee19}) and the uniform bounds obtained before, we get
\bFormula{conv1}
\Grad \Phi_\delta \to \Grad \Phi \ \mbox{in}\ L^q((0,T) \times \Omega; \RR^3) \ \mbox{for a certain}\ q > 2,
\eF
and
\bFormula{conv2}
\Phi_\delta \to \Phi \ \mbox{a.a. in} \ (0,T) \times \Omega.
\eF
Consequently, we can pass to the limit in (\ref{e1}), using the fact 
that $\Div \vu_\delta = S_{T,\delta}$ and the standard monotone 
operator theory to handle the limit in $\mu_\delta$.

\GGG Let us now test \eqref{e4} by $\phi\in W^{1,2}_0(\Omega)$. Then,
integrating by parts and using \eqref{ee17}, we easily arrive at
\bFormula{ee30}
 \int_0^T \| P_t \|^2_{W^{-1,2}(\Omega)} \ \dt \leq c.
\eF
Coupling this with \eqref{ee13}, we infer
\begin{equation}\label{conv21a}
  P_\delta \to P \ \mbox{strongly in}\ L^2(0,T:W^{-\epsilon,2}(\Omega)) \ 
   \text{for every }\epsilon\in (0,1).
\end{equation}
Let now $\phi\in C^\infty(\RR^3)$ with support in $\Omega$. Then, from 
\eqref{conv21a} and \eqref{ee20}, we obtain
\[
  \int_0^T \intO{ P_\delta \vu_\delta \phi } \ \dt
    \to \int_0^T \intO{ P \vu \phi } \ \dt,
\]
whence we can identify the limit of the product \EEE
\begin{equation}\label{conv21}
  \vu_\delta P_\delta \to \vu P \ \mbox{weakly in}\ L^2((0,T) \times \Omega; \RR^3).
\end{equation}
\GGG Next, testing \eqref{e5} by $n_\delta$ and using \eqref{ee13--} and \eqref{ee13},
it is easy to infer
\begin{equation}\label{conv21b}
  n_\delta \to n \ \mbox{weakly in}\ L^2(0,T;W^{1,2}_0(\Omega)),
\end{equation}
whence, using \eqref{conv21a} again, \EEE
\begin{align}\label{conv31}
  &P_\delta n_\delta \to P n \ \mbox{weakly-(*) in}\ L^\infty ((0,T) \times \Omega), \\
 \label{conv32}
  &\ P_\delta b(n_\delta) \to P \Ov{ b(n) } \ \mbox{weakly-(*) in} \ L^\infty ((0,T) \times \Omega),
\end{align}
for any \GGG $C^1$ \EEE function $b$, 
where $\Ov{b(n)}$ denotes a weak limit of $\{ b(n_\delta) \}_{\delta > 0}$.


\subsection{Strong convergence of the nutrients}

We finish the proof of compactness by showing strong (a.a.) 
pointwise convergence of the nutrients $\{ n_\delta \}_{\delta > 0}$. We have
\begin{equation}\label{conv99}
 -\Delta n_\delta + P n_\delta + \left[ \nu_1 (1 - Q(\Phi)) + \nu_2 Q(\Phi) \right] n_\delta = 
\end{equation}
\[
 = (P - P_\delta) n_\delta
  + \left[ \nu_1 (1 - Q(\Phi_\delta)) + \nu_2 Q(\Phi_\delta) \right] n_c
\]
\[
 + \Big( \left[ \nu_1 (1 - Q(\Phi)) + \nu_2 Q(\Phi) \right] - \left[ \nu_1 (1 - Q(\Phi_\delta)) + \nu_2 Q(\Phi_\delta ) \right] \Big)n_\delta,
\]
and, for the limit system,
\[
-\Delta n + P n + \left[ \nu_1 (1 - Q(\Phi)) + \nu_2 Q(\Phi) \right] n = \left[ \nu_1 (1 - Q(\Phi)) + \nu_2 Q(\Phi) \right] n_c.
\]
Thus, testing respectively by $n_\delta$ and $n$, integrating by parts,
and making use of the relations (\ref{conv1}--\ref{conv32}) 
\GGG (in particular, \eqref{conv32} is exploited with the choice $b(n_\delta)=n_\delta^2$ in order to 
manage the first term on the \rhs\ of \eqref{conv99}), \EEE 
we may show that
\[
\int_0^T \intO{ |\Grad n_\delta |^2 + P n^2_\delta + \left[ \nu_1 (1 - Q(\Phi)) + \nu_2 Q(\Phi) \right] n^2_\delta } \ \dt
\]
\[
\to \int_0^T \intO{ |\Grad n |^2 + P n^2 + \left[ \nu_1 (1 - Q(\Phi)) + \nu_2 Q(\Phi) \right] n^2 } \ \dt,
\]
which yields the desired conclusion
\bFormula{conv4}
\Grad n_\delta \to \Grad n, \ n_\delta \to n \ \mbox{in}\ L^2((0,T) \times \Omega).
\eF


\section{Approximation scheme}

\label{A}

In this section we  briefly introduce the approximated scheme  needed
to   obtain  rigorously the above described  a priori estimates. This part is quite standard, hence 
some details are omitted.


\subsection{Local existence by fixed point argument}

\newcommand{\ovu}{\Ov{\vu}}
\newcommand{\ovS}{\Ov{S}}
\newcommand{\oS}{\Ov{S}}
\newcommand{\oPhi}{\Ov{\Phi}}
\newcommand{\omu}{\Ov{\mu}}

Let 
\bFormula{regooS}
  \oS \in L^8(0,T;L^2(\Omega)), \
  \| \oS \|_{L^8(0,T;L^2(\Omega))}\le R_2,
\eF
where $R_2 > 0$ (this value can  be chosen
arbitrarily; for instance we can take $R_2=1$).

Replace $S_T +\lambda_3 \Phi$ with $\oS$ and 
solve (\ref{e1}--\ref{e3}) locally in
time by  a fixed point argument.
The following can be proven:
\bLemma{lemma1}
 Let $\oS$ be given by~(\ref{regooS}). Let  $\delta\in(0,1/4)$ and 
 $\Phi_{0\delta}\in W^{2,6}(\Omega),$ 
 \GGG $\Phi_{0\delta}\in [\delta,1-\delta]$. \EEE
 Then there exists $T_0\in(0,T]$ possibly depending on $\delta$
 such that the system
 \bFormula{e1d}
  {a.\;\;\;}\partial_t \Phi - \delta \Delta \mu_t + \vu \cdot \Grad \Phi - \Delta \mu = 0, \ \ \   
    {b.\;\;\;} \mu = - \Delta \Phi + \mathcal{F}'(\Phi),
 \eF
 \bFormula{e2d}
   \vu = - \Grad \Pi + \mu \Grad \Phi,
 \eF
 \bFormula{e3d}
  - \Delta \Pi = - \Div ( \mu \Grad \Phi ) + \oS - \lambda_3 \Phi,
 \eF
 coupled with the initial and boundary conditions
 \begin{equation}\label{boudelta}
 \mu = \Pi = 0, \quad \Grad\Phi\cdot\nu=0\quad\hbox{on }\partial\Omega\times (0,T),
 \end{equation}
 \begin{equation}\label{inidelta}
 \mu(0)=0, \quad \Phi(0)=\Phi_{0\delta},
 \end{equation}
 has at least one solution $(\Phi,\mu,\Pi,\vu)$ satisfying
 the regularity properties
 \begin{align}\label{rlemma11}
  &\Phi\in  H^1(0,T_0;H^1(\Omega))\cap L^\infty(0,T_0;W^{2,6}(\Omega)),\\
 \label{rlemma12}
   & \mu\in H^1(0,T_0;H^1(\Omega))\cap L^\infty(0,T_0;H^2(\Omega)),\\
 \label{rlemma13}
 & \Pi\in L^8(0,T;H^2(\Omega))\,. 
 \end{align}
\eL
\noindent%
\begin{proof}
Let $T_0\in(0,T]$ to be chosen below and let
\[
  \oPhi \in L^{4}(0,T_0;W^{1,4}(\Omega)), \quad
   \omu \in L^{4}((0,T_0)\times\Omega),
\]
with
\[
  \| \oPhi \|_{L^{4}(0,T_0;W^{1,4}(\Omega))}
   + \| \omu \|_{L^{4}((0,T_0)\times\Omega)} \le R_1.
\]
This in particular implies
\[
  \| \omu \Grad \oPhi \|_{L^2((0,T_0)\times\Omega)}
   \le Q(R_1).
\]
Again, $R_1>0$ can  be chosen  arbitrarily.
Here and below $Q$ is a computable positive function, monotone
increasing  in each of its arguments.

Replace $\oPhi$, $\omu$ and $\oS$ in equation
\eqref{e3d}, which  becomes
\bFormula{e2fp}
  - \Delta \Pi = - \Div ( \omu \Grad \oPhi ) + \oS - \lambda_3 \oPhi
\eF
and is still endowed with the boundary condition $\Pi=0$.
Clearly, \eqref{e2fp} has one and only one solution
\bFormula{rego:Pi}
  \Pi \in L^2(0,T_0;H^1_0(\Omega)).
\eF
Moreover,
\[
  \| \Pi \|_{L^2(0,T_0;H^1_0(\Omega))} \le Q(R_1,R_2).
\]
Set
\bFormula{e2b}
  \vu := - \Grad \Pi + \omu \Grad \oPhi \in L^2((0,T_0)\times \Omega; \RR^3),
\eF
and replace it in \eqref{e1d}. 
Once $\vu$ is assigned we can easily prove existence of a solution to
\eqref{e1d}. Note that no regularization of
$\mathcal{F}$ is required. The regularity class of the solution can be formally
determined multiplying {\eqref{e1d} a.} by $\mu$ and  {\eqref{e1d} b.} by $\Phi_t$.
Note that
\[
 \left| \intO{ \vu \cdot \Grad \Phi \mu } \right|
   \le \| \vu \|_{L^2(\Omega;\RR^3)} \| \mu \|_{L^4(\Omega)}
    \| \Grad \Phi \|_{L^4(\Omega;\RR^3)}
\]
\[
   \le c \| \vu \|_{L^2(\Omega;\RR^3)}^2 \| \mu \|_{H^1(\Omega)}^2
     + c \| \Phi \|_{H^2(\Omega)}^2
\]
\[
     \le c \| \vu \|_{L^2(\Omega;\RR^3)}^2 \| \mu \|_{H^1(\Omega)}^2
     + c \| \mu \|_{L^2(\Omega)}^2
     + c \| \Phi \|_{H^{1}(\Omega)}^2,
\]
The last inequality follows by multiplying  \eqref{e1d} {b.}
by $\Delta \Phi$ and using the monotonicity of $\mathcal{C}'$
(cf.~Subsec.~\ref{ass}).
Then, we can apply Gr\"onwall's lemma to obtain
\[
  \| \Phi \|_{L^\infty(0,T_0;H^1(\Omega))}
   + \| \mu \|_{L^\infty(0,T_0;H^1(\Omega))}
   \le Q(R_1,R_2,\delta^{-1},T_0).
\]
Next, multiplying  \eqref{e1d} {a.}\
by $\mu_t$, the time derivative of {\eqref{e1d} b.}\
by $\Phi_t$, and summing the results, yields
\bFormula{g11}
  \| \Phi \|_{H^1(0,T_0;H^1(\Omega))}
   + \| \mu \|_{H^1(0,T_0;H^1(\Omega))}
   \le Q(R_1,R_2,\delta^{-1},T_0).
\eF
Note that, to deduce \eqref{g11} in a rigorous way, it 
would have been necessary to regularize $\mathcal{C}$
in order for its second derivative to be well defined.
However this is a standard argument and the resulting estimate
would be independent of the regularization since it just relies on the
monotonicity of $\mathcal{C}'$. Hence, we omit giving details.
Finally, the same argument used for the complete system
yields  (cf.~\eqref{ee19})
\bFormula{regoPhifp}
  \| \Phi \|_{L^\infty(0,T_0;W^{2,6}(\Omega))}
   \le Q(R_1,R_2,\delta^{-1},T_0).
\eF
Next, multiplying {\eqref{e1d} a.}  by $\Delta \mu$ and using \eqref{regoPhifp},
it is not difficult to arrive at
\bFormula{regomufp}
  \| \mu \|_{L^\infty(0,T_0;H^{2}(\Omega))}
   \le Q(R_1,R_2,\delta^{-1},T_0).
\eF
By Sobolev's embedding this implies that there exists
$C_\delta>0$ such that
\bFormula{regomufp2}
  - C_\delta \le \mu(t,x) \le C_\delta
   \ \text{for a.e.~}(t,x) \in (0,T)\times \Omega.
\eF
\GGG Thanks to assumption \eqref{ee7c}, recalling that
$\Phi_{0\delta}\in [\delta,1-\delta]$, and 
applying maximum principle arguments in \eqref{e1d}, 
we deduce the following \emph{separation property}: \EEE
there exists $\kappa_\delta >0$ such that
\bFormula{separ}
  -1 + \kappa_\delta \le \Phi(t,x) \le 1- \kappa_\delta
   \ \text{for a.e.~}(t,x) \in (0,T)\times \Omega.
\eF
Thanks to the above separation property, we can
easily prove the uniqueness of the couple
$(\Phi,\mu)$. We already know that, once $\oS$ is given,
a unique $\Pi$ solving \eqref{e2fp} is determined. Hence
we have a unique $\vu$ given by \eqref{e2b}. Assuming that, for this
$\vu$, a couple of pairs $(\Phi,\mu)$ solve \eqref{e1d},
we can test the difference of {\eqref{e1d} a.} by the difference of the $\mu$'s
and the difference of {\eqref{e1d} b.}  by
the difference of the $\Phi_t$'s. Performing standard manipulations
and using the separation property \eqref{separ}
it is then easy to deduce a contractive estimate.
Hence, the couple $(\Phi,\mu)$ is in fact unique.

The above argument permits us to define, for the {\em fixed}\/
$\oS$ given by \eqref{regooS}, the map
\[
  \calM_1: B_{R_1} \to L^\infty(0,T_0;W^{2,6}(\Omega)) \times L^\infty(0,T_0;H^{2}(\Omega)),
   \ \calM_1:(\oPhi,\omu) \mapsto (\Phi,\mu),
\]
where $B_{R_1}$ is the closed ball of radius $R_1$ in the space
$L^{4}(0,T;W^{1,4}(\Omega)) \times L^{4}((0,T)\times\Omega)$.
We aim to apply Schauder's fixed point theorem to the above map.
To this purpose, we first observe that we can take $T_0$ small
enough so that the map takes values into $B_{R_1}$. Moreover, $\calM_1$
is compact by Sobolev's embeddings. Finally, the continuity
of $\calM_1$ can be shown by standard methods relying on the a priori
estimates obtained above.

Hence, by Schauder's theorem, there exists a time $T_0\le T$,
possibly depending on $\delta$, such that system (\ref{e1d}--\ref{e3d}),
coupled with the initial and boundary conditions, has at least one solution
$(\Phi,\mu,\Pi,\vu)$, 
in the interval $(0,T_0)$. The regularity of this solution is 
 specified by \eqref{rego:Pi},
\eqref{g11}, \eqref{regoPhifp}, and \eqref{regomufp}.
To conclude the proof, it  remains to improve the
regularity of $\vu$. Since 
we now know that $\Pi$ solves \eqref{e3d}, using
\eqref{regoPhifp} and \eqref{regomufp} it is easy to check that
\[
  \| - \Div (\mu \Grad\Phi) \|_{L^\infty(0,T_0;L^6(\Omega))}
   \le Q(R_2,\delta^{-1},T_0).
\]
Hence, recalling \eqref{regooS} and applying elliptic regularity
to \eqref{e3d}, we arrive at
\bFormula{regoufp}
  \| \Pi \|_{L^8(0,T_0;H^2(\Omega))} +
   \| \vu \|_{L^8(0,T_0;H^1(\Omega;\RR^3))}    \le Q(R_2,\delta^{-1},T_0).
\eF
This gives \eqref{rlemma13} and concludes the
proof of Lemma~\ref{Llemma1}.
\end{proof}
\GGG \bRemark{rem:sep}
 As already noted in Remark~\ref{RRsing}, assumption \eqref{ee7c} is 
 needed only for the sake of obtaining higher regularity 
 of approximating functions. Indeed, it would be enough to 
 assume it to hold in the approximation (for a suitable 
 family $\mathcal{C}_\delta$ tending to $\mathcal{C}$ as $\delta\to 0$)
 and not necessarily for $\mathcal{C}$.
\eR \EEE
\noindent%
\bLemma{lemma1b}
 For any $\oS$ as in (\ref{regooS}), the 
 quadruple $(\Phi,\mu,\Pi,\vu)$ solving (\ref{e1d}--\ref{e3d}) with the 
 initial and boundary conditions (\ref{boudelta}--\ref{inidelta}) 
 is~{\rm unique}.
\eL
\begin{proof}
A contractive estimate can be obtained simply by multiplying the
difference of {\eqref{e1d} a.} 
by the difference of the $\mu$'s,
the difference of {\eqref{e1d} b.} by
the difference of the $\Phi_t$'s, and the difference
of the \eqref{e3d} by the difference of the $\Pi$'s. We leave the details to the reader. We
note  that the separation property \eqref{separ}
and the additional regularity \eqref{regoufp} play a
role in this argument.
\end{proof}
\noindent%
Thanks to the above Lemmas, given $\oS$, there exists a {\em unique}\/
quadruple $(\Phi,\mu,\Pi,\vu)$ solving (\ref{e1d}--\ref{e3d}).
We now plug this quadruple into (a proper
regularization of) system (\ref{e4}--\ref{e5}). Namely,
we have the
\bLemma{lemma2}
 Let $\oS$ as in~(\ref{regooS}) and let $T_0$,
 $\Phi$, $\mu$, $\Pi$ and $\vu$ be given
 by Lemma~\ref{Llemma1}. Let $P_{0\delta}\in H^1_0(\Omega)$,
 $P_{0\delta}\in [0,1]$ a.e.\,. Then there exists one and
 only one couple $(P,n)$ satisfying the system
 \bFormula{e4fp}
   \partial_t P - \delta \Delta P + \Div (\vu P) = \Phi
    \big( n - \lambda_1 - \lambda_2 H (n_N - n) \big) P,
 \eF
 \bFormula{e5fp}
   - \Delta n + nP = T_c (n, \Phi),
 \eF
 over $(0,T_0)$, together with the initial and boundary conditions
 specified at the beginning (with $P_{0}$ replaced by $P_{0\delta}$)
 and the additional condition
 \bFormula{DirP}
   \delta P = 0 \ \mbox{on}\ \partial \Omega.
 \eF
 Moreover,
 \bFormula{maxpr}
   {a.} \ \ P(t,x) \ge 0,\,\,\, {b.} \ \ 0 \le n(t,x) \le 1 \
   \mbox{for a.e.}\ (t,x) \in (0,T_0) \times \Omega
 \eF
 and  the following regularity properties hold
 \bFormula{rego:Pfp}
  \| P \|_{H^1(0,T_0;L^2(\Omega))}
   + \| P \|_{L^\infty(0,T_0;L^2(\Omega))}
   + \| P \|_{L^2(0,T_0;H^2(\Omega))}
  \le Q(R_2,\delta^{-1},T_0),
 \eF
 \bFormula{rego:nfp}
  \| n \|_{H^1(0,T_0;H^1(\Omega))}
   + \| n \|_{L^\infty(0,T_0;H^2(\Omega))}
  \le Q(R_2,\delta^{-1},T_0),
 \eF
\eL
\noindent%
\begin{proof}
Let us introduce the truncation operator
$\calT(r)=\max\{0,\min\{1,r\}\}$. Plugging $\calT$ into the
\rhs\ of \eqref{e4fp}, we obtain the elliptic-parabolic system
\bFormula{e4tron}
  \partial_t P - \delta \Delta P + \Div (\vu P) = \Phi
   \big( \calT(n) - \lambda_1 - \lambda_2 H (n_N - n) \big) P,
\eF
\bFormula{e5tron}
  - \Delta n + nP = T_c (n, \Phi),
\eF
Existence of solutions to (the initial-boundary value problem) for
(\ref{e4tron}--\ref{e5tron}) is standard. For instance, one may
prove it by using the Faedo-Galerkin scheme.
Hence, we omit the details. Rather, we point
out which are the main a priori estimates involved, with
the  purpose of establishing sufficient regularity properties
of solutions. We will also see, as a byproduct, that the component
$n$ turns out to
take values in the interval $[0,1]$ so that the couple
$(n,P)$ will in fact solve (\ref{e4fp}--\ref{e5fp})
(without truncation).

To carry out this program, we start with multiplying
\eqref{e4tron} by $P$ to get
\[
  \frac12\ddt \| P \|_{L^2(\Omega)}^2
   + \delta \| \Grad P \|_{L^2(\Omega;\RR^3)}^2
\]
\[
   \le \intO{P \vu \cdot \Grad P}
    + \intO{\Phi
   \big( \calT(n) - \lambda_1 - \lambda_2 H (n_N - n) \big) P^2}.
\]
where we used  condition \eqref{DirP}.
In view of the smoothness of $\Phi$ and the
presence  of the truncation operator, the only
term that needs to be estimated is the first one on the \rhs.
By Poincar\'e's and Young's inequalities, we have
\begin{align}\label{gron}
  \intO{P \vu \cdot \Grad P}
   & \le \| P \|_{L^4(\Omega)} \| \vu \|_{L^4(\Omega;\RR^3)} \| \Grad P \|_{L^2(\Omega)}
\\
\nonumber
  & \le \| P \|_{H^1(\Omega)}^{7/4} \| P \|_{L^2(\Omega)}^{1/4} \| \vu \|_{L^4(\Omega;\RR^3)} \\
 \nonumber
   & \le \frac\delta2\| \Grad P \|_{L^2(\Omega)}^{2} + c_\delta \| P \|_{L^2(\Omega)}^{2} \| \vu \|_{L^4(\Omega;\RR^3)}^8.
\end{align}
Hence, by  \eqref{regoufp} and Gr\"onwall's Lemma, we arrive at
\bFormula{regoPfp}
  \| P \|_{L^\infty(0,T_0;L^2(\Omega))}
   + \| P \|_{L^2(0,T_0;H^1(\Omega))}
   \le Q(R_2,\delta^{-1},T_0).
\eF
Combining this relation with \eqref{regoufp} we infer 
\bFormula{Pfp2}
  \| \Div (\vu P) \|_{L^{8/5}(0,T_0;L^{3/2}(\Omega))}
   \le Q(R_2,\delta^{-1},T_0),
\eF
whence, applying parabolic regularity theory to \eqref{e4tron},
\bFormula{Pfp3}
  \| P_t \|_{L^{8/5}(0,T_0;L^{3/2}(\Omega))}
   + \| P \|_{L^{8/5}(0,T_0;W^{2,3/2}(\Omega))}
  \le Q(R_2,\delta^{-1},T_0).
\eF
In turn, by interpolation, this gives
\bFormula{Pfp4}
  \| P \|_{L^{8/3}(0,T_0;W^{\frac32-\epsilon,3/2}(\Omega))}
  \le Q(R_2,\delta^{-1},T_0,\epsilon^{-1}).
\eF
for all $\epsilon\in(0,1)$, whence, by Sobolev's embeddings,
\bFormula{Pfp5}
  \| P \|_{L^{8/3}(0,T_0;L^{6-\epsilon}(\Omega))}
   + \| \Grad P \|_{L^{8/3}(0,T_0;L^{2-\epsilon}(\Omega))}
  \le Q(R_2,\delta^{-1},T_0,\epsilon^{-1}).
\eF
Consequently, by \eqref{regoufp},
\bFormula{Pfp2b}
  \| \Div (\vu P) \|_{L^{2}(0,T_0;L^{\frac32-\epsilon}(\Omega))}
   \le Q(R_2,\delta^{-1},T_0,\epsilon^{-1}),
\eF
whence, going back to \eqref{e4tron},
\bFormula{Pfp3b}
  \| P_t \|_{L^{2}(0,T_0;L^{\frac32-\epsilon}(\Omega))}
   + \| P \|_{L^{2}(0,T_0;W^{2,\frac32-\epsilon}(\Omega))}
  \le Q(R_2,\delta^{-1},T_0,\epsilon^{-1}).
\eF
We  proceed by  bootstrapping. Actually,
some more iterations (whose details
are omitted for brevity) permit us to obtain \eqref{rego:Pfp}.
Once sufficient regularity is achieved, the same maximum
principle argument used for the coupled system
gives  \eqref{maxpr} {a.}

\smallskip

We  now pass to equation \eqref{e5tron}.
By elliptic regularity (i.e., multiplying
by $n-1\in H^1_0(\Omega)$), we infer
\bFormula{regonfp2}
  \| n \|_{L^{2}(0,T_0;H^1(\Omega))}
   \le Q(R_2,\delta^{-1},T_0).
\eF
Next, using, as for the complete system, the sign
condition on the \rhs, we get the second \eqref{maxpr}.
This entails in particular that $P$ solves
\eqref{e4fp}, i.e., no truncation in fact occurs.

\noindent By \eqref{regoPfp}, \eqref{maxpr}, and elliptic regularity,
it follows that
\bFormula{regonfp3}
  \| n \|_{L^{\infty}(0,T_0;H^2(\Omega))}
   \le Q(R_2,\delta^{-1},T_0).
\eF
Now, we differentiate \eqref{e5fp} in time. Recalling
\eqref{def:Tc} we  have
\bFormula{e5fpt}
  - \Delta n_t + n_t P + n P_t
   = - \left[ \nu_1 (1 - Q(\Phi)) + \nu_2 Q(\Phi) \right] n_t
\eF
\[
  \mbox{} + \left[ - \nu_1 Q'(\Phi) + \nu_2 Q'(\Phi) \right] \Phi_t (n_c - n).
\]
Test the above relation by $n_t$ and use the regularity given by
\eqref{g11} and \eqref{rego:Pfp} together with the second~\eqref{maxpr} and 
the positivity of the given term $\left[ \nu_1 (1 - Q(\Phi)) + \nu_2 Q(\Phi) \right]$,
to obtain 
\bFormula{regonfp4}
  \| n_t \|_{L^{2}(0,T_0;H^1(\Omega))}
   \le Q(R_2,\delta^{-1},T_0).
\eF
This, combined with \eqref{regonfp3}, yields \eqref{rego:nfp}.
Finally, we have to prove uniqueness of the solution $(n,P)$.
This in fact follows from a standard argument. Indeed, it is sufficient
to test the difference of \eqref{e4fp} by the difference of the
$P$'s and the difference of \eqref{e5fp} by the difference of the
$n$. Then, the transport term in \eqref{e4fp} is treated in a way similar
to~\eqref{gron}, whereas the \rhs\ of \eqref{e5fp} is easily
controlled in view of the high regularity of $\Phi$ and of the
sign condition. This concludes the proof of the lemma.
\end{proof}

\smallskip

\noindent%
We can now finalize our fixed point argument for the complete system.
\bTheorem{theorem1}
 Let $\delta\in (0,1/4)$, $\Phi_{0\delta}\in W^{2,6}(\Omega)$, $\Phi_{0\delta}\in [\delta,1-\delta]$,
 $P_{0\delta}\in H^1_0(\Omega)$, $P_{0\delta}\in[0,1]$ a.e.\,. 
 Then there exists $T_1\in(0,T]$
 possibly depending on $\delta$ such that the system
 \bFormula{e1r}
   \partial_t \Phi - \delta \Delta \mu_t + \vu \cdot \Grad \Phi - \Delta \mu = 0, \ \ \mu = - \Delta \Phi + \mathcal{F}'(\Phi),
 \eF
 \bFormula{e2r}
   \vu = - \Grad \Pi + \mu \Grad \Phi,
 \eF
 \bFormula{e3r}
  - \Div \vu =  S_T = nP - \lambda_3 (\Phi - P),
 \eF
 \bFormula{e4r}
   \partial_t P - \delta \Delta P + \Div (\vu P) = \Phi
    \big( n - \lambda_1 - \lambda_2 H (n_N - n) \big) P,
 \eF
 \bFormula{e5r}
   - \Delta n + nP = T_c (n, \Phi),
 \eF
 coupled with the initial and boundary conditions (\ref{e6}--\ref{e8})
 (with $\Phi_{0\delta}$ and $P_{0\delta}$ replacing $\Phi_{0}$ and $P_{0}$)
 and (\ref{DirP}), has at least one solution $(\Phi,\mu,\vu,P,n)$
 defined over the time interval $(0,T_1)$ and
 satisfying the regularity properties  (\ref{rlemma11}--\ref{rlemma12}), (\ref{regoufp}),
 (\ref{rego:Pfp}--\ref{rego:nfp}).
\eT
\begin{proof}
We let $\oS$ be as in \eqref{regooS}, where the choice of $R_2\ge 0$
is in fact arbitrary. Then, applying first
Lemmas~\ref{Llemma1}, \ref{Llemma1b}   and then Lemma~\ref{Llemma2} we obtain
a {\em unique}\/ quintuple $(\Phi,\mu,\vu,P,n)$.
Thus, we can consider the map
\bFormula{defiM2}
  \calM_2 : \oS \mapsto S := nP + \lambda_3 P.
\eF
In view of \eqref{rego:Pfp}, \eqref{rego:nfp} and Sobolev's
embeddings, it is easy to check that
\bFormula{rego:S}
  \| S \|_{L^\infty(0,T_0;L^2(\Omega))}
  \le Q(R_2,\delta^{-1},T_0),
\eF
In particular, we can choose $T_1 \in (0,T_0]$ such that
$S$ lies in the closed ball $B_{R_2}$ of $L^8(0,T_1;L^2(\Omega))$.
Moreover, continuity and compactness of the map
$\calM_2$ in the topology of $L^8(0,T_1;L^2(\Omega))$
are an easy consequence of the regularity properties
\eqref{rego:Pfp}, \eqref{rego:nfp}, the Lions-Aubin theorem,
and  the a priori estimates in
Lemmas~\ref{Llemma1}, \ref{Llemma2}. Hence we can apply once more
Schauder's theorem to $\calM_2$, which gives that
\bFormula{SoS}
  \oS = S = nP + \lambda_3 P
   \ \text{in}\ (0,T) \times \Omega.
\eF
Hence, (\ref{e2d}--\ref{e3d}) reduce to (\ref{e2}--\ref{e3}),
where $S_T$ is given by~\eqref{def:ST}. This concludes the proof
of the theorem.
\end{proof}
\noindent%
In order to complete the proof of Theorem~\ref{Ttheoremmain}, we need now 
to pass to the limit in the regularized system as $\delta\searrow 0$,
assuming of course that $\Phi_{0\delta} \to \Phi_0$ and 
$P_{0\delta} \to P_0$ in suitable ways. We just briefly comment 
on the most delicate part of this step, which consists in the passage to the limit in 
(\ref{e4r}) in order to recover (\ref{wf3}). The other parts are indeed standard since 
it can be immediately seen that the a-priori estimates performed in Section~\ref{B} 
are still valid on the regularized system and they turn out to be also independent of $\delta$.

Taking a test function $\varphi$ as in (\ref{wf3}) we multiply (\ref{e4r}) by $\varphi$ to obtain 
\bFormula{ident}
\int_0^T \intO{ \left[ P \partial_t \varphi + P \vu \cdot \Grad \varphi + 
\Phi \left(n - \lambda_1 - \lambda_2 H(n_N - n) \right) \varphi \right] } \ \dt 
\eF
\[
= \delta \int_0^T \intO{ \Grad P \cdot \Grad \varphi } \ \dt
- \intO{ P_0 \varphi (0, \cdot) } + \int_0^T \int_{\partial \Omega} P \vu \cdot \vc{n} \varphi \ {\rm dS}_x
\]
\[
- \delta \int_0^T \int_{\partial \Omega} \Grad P \cdot \vc{n} \varphi \ {\rm dS}_x, 
\]
where, as 
\[
P|_{\partial \Omega} = 0,\ P \geq 0 \ \mbox{in} \ (0,T) \times \Omega, 
\]
\[
\int_0^T \int_{\partial \Omega} P \vu \cdot \vc{n} \varphi \ {\rm dS}_x = 0,\ - \delta \int_0^T \int_{\partial \Omega} \Grad P \cdot \vc{n} \varphi \ {\rm dS}_x \geq 0.
\]
Letting $\delta \searrow 0$ in (\ref{ident}) we get (\ref{wf3}). 

Finally, let us notice that, by standard arguments, it is possible to show that 
the a priori estimates provide an extension of the local approximate solution 
up to the original final time $T$. Hence, in particular we have a global solution 
in the limit. This concludes the proof of Theorem~\ref{Ttheoremmain}.


\section{Singular limit}

\label{S}

In this section, we consider the problem  obtained from (\ref{e1}--\ref{e5}) by taking $S_T=S_D=0$. Hence we 
just consider the system of equations for $\Phi$ and $\vu$, decoupled from the rest, of the form
\bFormula{S1}
\partial_t \Phi + \Div(\vu  \Phi) - \Div (\Grad \mu ) =  0 , \ \mu = - \ep^2 \Delta \Phi + \mathcal{F}'(\Phi),
\eF
\bFormula{S2}
\vu = - \Grad \Pi + \mu \Grad \Phi,
\eF
\bFormula{S3}
\Div \vu = 0,
\eF
with the boundary conditions
\bFormula{S4}
\vu \cdot \vc{n}|_{\partial \Omega} = 0, \ \Grad \Phi \cdot \vc{n}|_{\partial \Omega} = 0, \ \mu|_{\partial \Omega} = 0.
\eF
\GGG Notice that, in particular, we are considering here a no-flux condition for $\Pi$ 
in place of the Dirichlet condition in~\eqref{e6}. \EEE

Similarly to Section \ref{ECH}, we derive the energy balance
\bFormula{S5}
\frac{{\rm d}}{{\rm d}t} \intO{ \left[ \frac{\ep^2}{2} |\Grad \Phi |^2 + \mathcal{F}(\Phi) \right] } +
\intO{ |\Grad \mu|^2 + |\vu|^2  } = 0.
\eF
Next,
\[
 \intO{ \left[ \ep^2 |\Delta \Phi |^2 + \mathcal{F}'' (\Phi) |\Grad \Phi |^2 \right] } = \intO{ \Grad \mu \cdot \Grad \Phi }.
\]
Then, assuming strict convexity of $\mathcal{F}$, namely
\bFormula{S6}
\mathcal{F}'' \geq \lambda > 0,
\eF
 the following estimates can be deduced
\bFormula{S7}
\int_0^T \| \ep \Delta \Phi \|^2_{L^2(\Omega)} \ \dt \leq c , 
 \ \int_0^T \| \Grad \Phi \|^2_{L^2(\Omega;\RR^3)} \ \dt \leq c.
\eF


\subsection{Compactness of the velocity}

In view of (\ref{S5}) we may assume there is a subsequence such that
\[
\vu_\ep\to \vu \ \mbox{weakly in}\ L^2((0,T) \times \Omega; \RR^3).
\]
Obviously,
\bFormula{S8}
  \Div \vu = 0, \ \vu \cdot \vc{n}|_{\partial \Omega} = 0.
\eF
We can now write
\[
\vue= - \Grad \left( \Pi_\ep - \mathcal{F}(\Phi_\ep) \right) - \ep^2 \Delta \Phi_\ep \Grad \Phi_\ep;
\]
whence, seeing that
\[
\ep^2 \Delta \Phi_\ep \Grad \Phi_\ep \to 0 \ \mbox{in} \ L^1 ((0,T) \times \Omega),
\]
we conclude that
\[
{\bf curl}_x \vu = 0,
\]
which, combined with (\ref{S8}), yields
\[
\vu = 0.
\]
Therefore, taking $\ep \to 0$, system (\ref{S1})--(\ref{S3}) converges to 
\bFormula{S1lim}
 \partial_t\Phi-\Delta\mu=0,  \qquad \mu=\mathcal F'(\Phi),
\eF
which satisfies the energy law
\bFormula{enlim}
  \ddt \intO {\mathcal F(\Phi)}
   +\intO {|\nabla_x\mu|^2} = 0.
\eF
\GGG Summarizing, we have proved the 
\bTheorem{theosing}
Let the assumptions given in Subsec.~\ref{ass} hold,
let $\mathcal{F}$ satisfy \eqref{S6}, and let 
$(\Phi_\ep,\mu_\ep,\vu_\ep)$ denote a family of weak solutions
to the system (\ref{S1}--\ref{S3}) complemented 
with the boundary conditions \eqref{S4}
and the Cauchy conditions. Then, as $\ep\to 0$, the functions 
$(\Phi_\ep,\mu_\ep,\vu_\ep)$ suitably tend to a triple 
$(\Phi,\mu,0)$ satisfying \eqref{S1lim} together with 
the energy equality \eqref{enlim}
and the initial and boundary conditions.
\eT
\bRemark{R1}
It would be interesting to investigate whether similar estimates could be derived
for the singular flux
\[
  \vu = - \Grad \Pi + \frac{1}{\ep} \mu \Grad \Phi.
\]
However, the above argument does not seem to be easily adaptable to cover such a situation.
For instance, we cannot prove uniform integrability of the product
\[
  \ep \Delta \Phi \Grad \phi
\]
in that case.
\eR
\EEE


\section{Acknowledgements}

The research of E.F.\ leading to these results has received funding 
from the European Research Council under the European Union's Seventh Framework Programme (FP7/2007-2013)/ ERC 
Grant Agreement 320078. The Institute of Mathematics of the Academy of Sciences of the Czech
Republic is supported by RVO:67985840. The work of E.R.\ and of G.S.\ was
supported by the FP7-IDEAS-ERC-StG Grant \#256872 (EntroPhase), by GNAMPA (Gruppo Nazionale per l'Analisi Matematica, 
la Probabilit\`a e le loro Applicazioni) of INdAM (Istituto Nazionale di Alta Matematica), and by IMATI -- C.N.R. Pavia.


\end{document}